
\magnification=1200 
\input amssym.def
\overfullrule0pt

\def\CC{{\cal C}}

\def\R{{\Bbb R}}
\def\D{{\Bbb D}}

\def\C{{\Bbb C}}

\def\11{{\bf 1\!\!}}
\def\ratop{\mathop{\hbox to .25in{\rightarrowfill}}\limits}

\def\nhang{\hangindent=4pc\hangafter=1}

\centerline{\bf  Upper semi-continuity of the Royden-Kobayashi
pseudo-norm,}

\centerline{\bf a counterexample for H\" olderian  almost complex structures.}

\bigskip
\bigskip

\centerline{
by  Sergey Ivashkovich, Sergey Pinchuk and Jean-Pierre Rosay%
\footnote*{Partly supported by NSF grant.}}

\footnote{} {AMS classification: 32Q60, 32Q65, 32Q45, 35J60.}
\bigskip\bigskip
If $X$ is an almost complex manifold, with an almost complex
structure $J$ of class $\CC^\alpha$, for some $\alpha >0$,
for every point $p\in X$ and every tangent vector $V$ at $p$,
there exists a germ of $J$-holomorphic disc through $p$ with
this prescribed tangent vector. This existence result goes back
to [N-W] (Theorem III.). See [I-R] (Appendix 1) for a re-writing and a
generalization (to $k$-jets) of the proof.
All the $J$ holomorphic curves are of class $\CC^{1,\alpha}$
([N-W], [S]).
\bigskip
Then, exactly as for complex manifolds one can define the Royden-Kobayashi
pseudo-norm of tangent vectors. The question arises whether this
pseudo-norm is an upper semi-continuous function on the tangent bundle.
For complex manifolds it is the crucial point in Royden's proof of the
equivalence of the two standard definitions of the Kobayashi pseudo-metric
([R 1], [R 2], [L] pages 88-94). The upper semi-continuity of the
Royden-Kobayashi pseudo-norm has been established by Kruglikov [K] for
structures that are smooth enough. In [I-R], it is shown
that $\CC^{1,\alpha}$ regularity of $J$ is enough.
\bigskip
Here we show the following:
\bigskip\noindent
{\bf Theorem 1.} {\it There exists an
almost complex structure $J$ of class
 $\CC^{1\over 2}$ on the unit bidisc $\D^2\subset \C^2$, such that the
Royden-Kobayashi pseudo-norm is not an upper semi-continuous
function on the tangent bundle.}
\bigskip
The example is very explicit and very simple to describe. See Part
II. We refer the reader unfamiliar with the above notions to
[I-R]. For the proof of the failure of upper semi-continuity we
shall need the following result:
\bigskip\noindent
{\bf Theorem 2.} {\it For any continuous (complex valued) function
$f$ defined on the unit disc in $\C$, that is continuously
differentiable on the set on which $f\neq 0$, and that on that set
satisfies: $${\partial f\over \partial \overline z}~=~|f|^{1\over
2}~, $$
one has ${\rm sup}_{|z|<1}|f(z)|\geq {1\over 10}$, if
$f(0)\neq 0$. }
\bigskip\noindent
So, roughly speaking, the theorem says that the equation
${\partial f\over \partial \overline z}~=~|f|^{1\over 2}~,$
which has of course  $f\equiv 0$ as a solution does not
have small solutions with $f(0)\neq 0$. It can of course
be compared with what happens for the ordinary differential equation
$g'=|g|^{1\over 2}$ on the interval $[-1~,~+1]$. (Up to a factor $2$, it
is the equation to which ${\partial f\over \partial \overline z}~=~|f|^{1\over 2}$
reduces if $f(z)$ depends only on the real part of $z$.)
Every (real) solution $g$ with $g(0)\neq 0$ satisfies $g(1)> {1\over 4}$ (if $g(0)>0$),
or $g(-1)<-{1\over 4}$ (if $g(0)<0$). However, the equation has small non zero solutions
that are identically 0 on large sub-intervals of $[-1~,~+1]$.
\bigskip
\centerline{\bf Part I.}
\bigskip
The crucial step for proving Theorem 2 is the following Lemma.
\bigskip
\noindent{\bf Lemma 1.} {\it Let $\omega$ be an open subset in $\C$.
Let $h$ be a continuously differentiable function defined on $\omega$
such that ${\partial h\over \partial \overline z} = |h|^{1\over 2}$.
If $h$ does not vanish anywhere on $\omega$, then $h$ is smooth
($\CC^\infty$) and
$$\Delta \big( |h|^{3\over 4})\geq {3\over 4} |h|^{-{1\over 4}}.$$}
\bigskip
\noindent {\bf Proof.} In the proof differentiation will be denoted by
lower indices.

The smoothness of $u$ follows from standard elliptic bootstrapping.
Since the question is purely local, we can assume that
$\omega$ is simply connected. Set $g=h^{1\over 2}$ (a determination
of the square root of $h$).
The hypothesis $h_{\overline z}=|h|^{1\over 2}$,
is equivalent to $\big( h^{1\over 2}\big)_{\overline z}=
{1\over 2}\big({\overline h\over h}\big)^{1\over 4}. $
It can then be restated:
$$g_{\overline z}={1\over 2}\big( {\overline g \over g}\big)^{1\over 2}. \eqno (1)$$
Set $g(z)=\rho e^{i\varphi (z)}$, with $\rho >0$ and $\varphi$ real valued.
The conclusion to be reached is
$$ \Delta (\rho^{3\over 2})\geq {3\over 4}\rho^{-{1\over 2}}. \eqno (*)$$
\bigskip
One has $g_{\overline z}=\rho_{\overline z}e^{i\varphi}+i\varphi_{\overline z} \rho
e^{i\varphi}$, $\big({\overline g\over g}\big)^{{1\over 2}}=e^{-i\varphi}$.
Hence (1) gives
$$\rho_{\overline z}+i\varphi_{\overline z}\rho={1\over 2}e^{-2i\varphi}.\eqno (2)$$
Separating real and imaginary parts:
$$\rho_x-\rho\varphi_y=\cos 2\varphi , \eqno (3)$$
$$\rho_y+\rho\varphi_x=-\sin 2 \varphi . \eqno (3')$$
Multiply the first equation by $-\varphi_y$, the second one by
$\varphi_x$ and add up in order to get:
$$\rho(\varphi_x^2+\varphi_y^2)+(\rho_y\varphi_x-\rho_x\varphi_y)
=-(\varphi_x\sin 2\varphi +\varphi_y \cos 2 \varphi).\eqno (4)$$
Now, differentiate (3) with respect to $x$ and (3') with respect
to $y$ and add up. One gets:
$$\Delta \rho+(\rho_y\varphi_x-\rho_x\varphi_y)
=-2(\varphi_x\sin 2\varphi + \varphi_y\cos 2\varphi).\eqno (5)$$
(4) and (5) yield:
$$\Delta \rho - 2\rho (\varphi_x^2+\varphi_y^2)=\rho_y\varphi_x -\rho_x\varphi_y,$$
and therefore
$$2\rho \Delta \rho - 4\rho^2 (\varphi_x^2+\varphi_y^2)=2\rho (\rho_y\varphi_x -\rho_x\varphi_y)
.\eqno (6)$$
From equation (3) we have:
$(\rho_x-\rho\varphi_y)^2+ (\rho_y+\rho\varphi_x)^2=1$, hence:
$$\rho_x^2+\rho_y^2+\rho^2(\varphi_x^2+\varphi_y^2)+2\rho
(\rho_y\varphi_x-\rho_x\varphi_y)=1.\eqno (7)$$
From (6) and (7):
$$\rho_x^2+\rho_y^2+2\rho\Delta \rho=1+3 \rho^2(\varphi_x^2+\varphi_y^2)\geq 1.\eqno (8)$$
The left hand side in (8) is equal to ${4\over 3}\rho^{1\over 2}
\Delta \big( \rho^{3\over 2} \big)$. This establishes (*).
\smallskip\hfill Q.E.D.
\bigskip
\noindent {\bf Lemma 2.} {\it  Let $u$ be a non-negative
subharmonic function on the unit disc in $\C$. Assume that on the
set on which $u\neq 0$, one has $\Delta u \geq 1$. If $u(0)>0$,
then ${\rm sup}_{|z|<1}u(z) > {1\over 4}$.}
\bigskip\noindent
{\bf Proof.} Set $\omega=\{ z\in \D;~u(z)\neq 0\}$. The function
$v=u-{1\over4}(x^2+y^2)$ is a subharmonic function on $\omega$.
If  ${\rm Sup}_{|z|<1}u(z) \leq {1\over 4}$, on the whole boundary of
$\omega$ we would have $v\leq 0$, which is impossible by the maximum principle
since $v(0)=u(0)>0$.
\smallskip\hfill Q.E.D.
\bigskip
\noindent {\bf Proof of Theorem 2.} Assume that ${\rm sup}~|f|<{1\over 10}$.
Then Lemma 1 gives $\Delta \big( |f|^{3\over 4}\big)\geq
{3\over 4} 10^{1\over4}\geq 1$ on $\omega$. So we can apply Lemma 2,
and one therefore gets $\sup_{|z|<1}|f(z)|^{3\over 4}\geq {1\over 4}$, and therefore
$\sup_{|z|<1}|f(z)|\geq {1\over 10}$, a contradiction, as desired.
\smallskip\hfill
Q.E.D.
\bigskip\bigskip

\centerline{\bf Part II.}
\bigskip
For $r>0$, $\D_r$ will denote the open disc $\{z=(x+iy)\in
\C;~|z|<r\}$. For the unit disc $\D_1$ in $\C$, we will abbreviate
to $\D$. We will identify $\R^4$, with coordinates
$(x_1,y_1,x_2,y_2)$ with $\C^2$ (setting $z_j=(x_j+iy_j)$). Let
$\Omega=\D_2\times \D_{1\over 10}$. On $\Omega$ consider the
almost complex structure $J$ defined by:

 $J = J(x_1,y_1,x_2,y_2) =    \left(\matrix{ 0  &   -1  &    0 & 0    \cr
                         1  &   0  &   0 & 0     \cr
                         0 &   \lambda  &    0 & -1     \cr
                          \lambda  &    0 & 1 & 0    \cr}\right)$
\medskip\noindent
where $\lambda (x_1,y_1,x_2,y_2)= -2(x_2^2+y_2^2)^{1\over 4}=-2|z_2|^{1\over 2}$.
Note that $J^2=-{\bf 1}$.
\bigskip
In the next lines, we use complex notations.
The map from $\D_2$ into  $\Omega$: $z\mapsto (z,0)$ is a $J$-holomorphic map.
This shows that the Royden-Kobayashi pseudo-norm of the  vector $(1,0)$
tangent at the point $(0,0)$, is $\leq {1\over 2}$.
\bigskip\noindent
CLAIM:{\it There exists $A<2$ such that if $Z$  is a $J$-holomorphic
map from $\D_r$ into $\Omega$ with:
$$Z(0)=(0,b)~{\rm with}~b\neq 0~,~{\rm and}~{\partial Z\over \partial x} (0)=(1,0)~,$$
then $r\leq A$.}
\bigskip
The claim implies that the vector $(1,0)$ tangent at the point
$(0,b)$ has Royden-Kobayashi pseudo-norm $\geq {1\over A} >
{1\over 2} $. So, by identification of $\Omega$ and $\D^2$,
to prove Theorem 1, it is enough to check the
claim.
\bigskip
\noindent
{\bf Proof of the Claim.}  To say that a map $z\mapsto Z(z)=(Z_1(z),Z_2(z))$, from
$\D_r$ to $\Omega$ is $J$-holomorphic, means that
$${\partial Z\over \partial y}(z)~=~
J(Z(z))~{\partial Z\over \partial x}(z).$$
It immediately gives us that $Z_1$ must be a holomorphic function of
$z$. By the Schwarz Lemma, if $r$ is close enough to 2, and $Z_1(0)=0$,
$Z_1'(0)=1$, $Z_1$ must be close to the identity, say on the disc of radius
${3\over 2}$. More precisely, there exists $A<2$ so that if $r>A$
there is a holomorphic map $\varphi$ defined on the unit disc, close
to the identity, with $\varphi (0)=0$ and $\varphi'(0)=1$ such that
$Z_1(\varphi (z))\equiv z$. Then consider the $J$-holomorphic map from
$\D$ into $\Omega$ obtained by re-parameterization:
$$z\mapsto Z^\#(z)=Z(\varphi (z))=(z,f(z)).$$
In order to prove the claim, we need to show that such a $J$-holomorphic map
$Z^\#$ does not exist with $f(0)\neq 0$.
\bigskip
We now look at the condition ${\partial Z^\#\over \partial y}(z)
=J(Z^\# (z)) {\partial Z^\#\over \partial x}(z)$, which is the condition
for $Z^\#$ to be $J$-holomorphic. Set $f=u+iv$.
Looking at the $(x_2,y_2)$ coordinates, one gets:
$${\partial u\over \partial y} = -{\partial v\over \partial x},$$
$${\partial v\over \partial y} = \lambda + {\partial u \over \partial x}.$$
Multiply the first equation by $i$ and subtract from the second equation.
One gets ${\partial (u+iv)\over \partial \overline z}=-{\lambda\over 2}$.
\smallskip
So the condition for the map $Z^\#$ be $J$-holomorphic, with $J$ as above,
is simply that ${\partial f \over \partial \overline z}=|f(z)|^{1\over 2}$.
If $f(0)\neq 0$, Theorem 2 implies that
${\rm Sup}_{|z|<1}|f(z)| > {1\over 10}$, which is impossible since
$Z^\#$ maps $\D$ into $\Omega$.
\medskip\hfill Q.E.D.
\bigskip
\noindent {\bf Comments.} \hfill\break
1) It is well known that for rough (non $\CC^1$)
almost complex structures
there is no uniqueness result: Two $J$-holomorphic maps defined on
$\D$ can coincide on a non empty open set without being identical.
It indeed happens in the above example. With the above notations,
taking $f$ depending only on $x$, it reduces to the non-uniqueness property
for the O.D.E. $g'=|g|^{1\over 2}$.
An example of $J$-holomorphic map $Z^\#$ (defined near 0) is then given by taking
$f(z)=0$ if $x\leq 0$, and $f(z)= x^2$ for $x>0$. 
A strong relation seems to exist between the lack of uniqueness for O.D.E. and the lack
of `small solutions'.
\bigskip\noindent
2) A minor modification of the proof leads to an example with an almost complex
structure H\" olderian of H\"older exponent ${2\over 3}$, instead of only ${1\over 2}$.
\bigskip
\noindent 3) Interestingly, there is another situation somewhat comparable
to the situation in this paper, and to which
L. Lempert drew the attention of one of us a few years ago. In [C-R] there is a
very simple construction of holomorphic discs, depending on a parameter,
simply using the implicit function Theorem, in complex dimension 2.
L. Lempert pointed out that, in higher dimensions, one could still prove the
existence of similar discs, by applying the Schauder fixed point theorem.
But it did not give the dependence on parameters as
needed in Kontiuit\" atsatz arguments. It could not due to the failure of the
Hartogs-Chirka Theorem in complex dimension $>2$ ([Ro]).
For $J$-holomorphic curves, proving the existence of curves with prescribed tangent
can be done with the implicit function Theorem if $J$ is of class $\CC^{1,\alpha}$.
But for $\CC^\alpha$ regularity of $J$, the proof as in [N-W] uses the Schauder fixed
point Theorem. The non upper semi-continuity of the Royden-Kobayashi pseudo-norm
comes from the impossibility of small perturbations.
\bigskip
\centerline{REFERENCES.}
\bigskip
\item{[C-R]} E.\ Chirka, J-P.\ Rosay. Remarks on the proof of a generalized
Hartogs Lemma. Ann. Pol. Math. {\bf LXX} (1998), 43-47.
\smallskip

\item{[I-R]} S.\ Ivashkovich, J-P.\ Rosay. Schwarz-type Lemmas for solutions
of $\overline \partial$-inequalities and complete hyperbolicity
of almost complex manifolds. (Preprint)
\smallskip
\item{[K]} B.S.\ Kruglikov. Existence of Close Pseudoholomorphic
Disks for Almost Complex Manifolds and an Application to the
Kobayashi-Royden Pseudonorm.
Funct. Anal. and  Appl. {\bf 33}
(1999) 38--48.
\smallskip
\item{[K-O]} B.S.\ Kruglikov, M. Overholt. Pseudoholomorphic mappings
and Kobayashi hyperbolicity. Differential Geom. Appl. {\bf 11}
(1999), 265--277.
\smallskip
\item{[L]} S.\ Lang. {\it Introduction to Complex Hyperbolic Spaces}
Springer-Verlag (1987).
\smallskip
\item{[N-W]} A.\ Nijenhuis, W.\ Woolf. Some integration problems in
almost complex and complex manifolds. Ann. of Math. {\bf 77} (1963), 424--489.
\smallskip
\item{[Ro]} J-P.\ Rosay. A counterexample to the Hartogs phenomenon
(a question by E. Chirka). Michigan Math. J. {\bf 45} (1998),
529-535.
\smallskip

\item{[R 1]} H.\ Royden. Remarks on the Kobayashi metric. Proc.
Maryland Conf. on Several Complex Variables. Springer Lecture Notes in
Math. {\bf 185}
(1971), 125-137.
\smallskip
\item{[R 2]} H.\ Royden. The extension of regular holomorphic maps.
Proc. A.M.S. {\bf 43} (1974), 306-310.
\smallskip

\item{[S]} J.-C. Sikorav. Some properties of holomorphic curves in
almost complex manifolds. In {\it Holomorphic Curves in Symplectic
Geometry}, eds. M. Audin and J. Lafontaine, Birkhauser (1994), 351--361.
\smallskip

\bigskip

\nhang{S.\ Ivashkovich: D\'epartement de Math\'ematiques, Universit\'e Lille I,
59655 Villeneuve d' Asq Cedex, France. {\it ivachkov@gat.univ-lille1.fr}}
\smallskip

\nhang{S.\ Pinchuk: Department of Mathematics, Indiana University,
Bloomington IN 47405 USA. {\it pinchuk@indiana.edu}}
\smallskip

\nhang{J-P.\ Rosay: Department of Mathematics, University of Wisconsin,
Madison WI 53706 USA. {\it jrosay@math.wisc.edu}}
\bye